\documentclass[twoside,10pt,reqno]{amsart}
\usepackage{amsmath,amssymb,amscd,verbatim}

\newcommand{\Z}{\ensuremath{\mathbb{Z}}}
\newcommand{\U}{\ensuremath{U(\mathfrak{g})}}
\newcommand{\Uh}{\ensuremath{U_{h}(\mathfrak{g}) } }

\newcommand{\Uhb}{\ensuremath{U_{h}(\mathfrak{b_{+}}) } }
\newcommand{\Uhbhat}{\ensuremath{U_{h}(\tilde{\mathfrak{b}}_{+}) } }
\newcommand{\UhA}{\ensuremath{U_{h}(\mathfrak{g},A,\tau) } }
\newcommand{\UhAp}{\ensuremath{U_{h}(\mathfrak{g}',A',\tau') } }

\newcommand{\Hom}{\ensuremath{\operatorname{Hom}}}

\newcommand{\Aut}{\ensuremath{\operatorname{Aut}}}
\newcommand{\g}{\ensuremath{\mathfrak{g}}}
\newcommand{\ghat}{\ensuremath{\tilde{\mathfrak{g}}}}

\newcommand{\h}{\ensuremath{\mathfrak{h}}}
\newcommand{\Uhh}{\ensuremath{U_{h}(\mathfrak{\h}) } }
\newcommand{\nilp}{\ensuremath{\mathfrak{n}}}
\newcommand{\borel}{\ensuremath{\mathfrak{b}}}
\newcommand{\borelp}{\ensuremath{\mathfrak{b}_{+}}}
\newcommand{\borelhat}{\ensuremath{\tilde{ \mathfrak{b}}}}

\newcommand{\C}{\ensuremath{\mathbb{C}} }

\newcommand{\p}[1]{\ensuremath{\bar {#1}}}

\newcommand{\D}[1]{\ensuremath{\Delta({#1})} }

\newcommand{\Uhatp}{\mathcal{\widetilde{U}_{+}}}
\newcommand{\Up}{\mathcal{U_{+}}}
\newcommand{\Upqd}{\mathcal{U_{+}^{\star}}}
\newcommand{\cas}{\Omega }

\newcommand{	\qd}{{}^{\star}}
\newcommand{	\qdo}{{}^{\star op}}

\newcommand{\DJgA}{\ensuremath{U_{h}^{DJ}(\g,A,\tau) } }
\newcommand{\DJgAp}{\ensuremath{U_{h}^{DJ}(\g',A',\tau') } }
\newcommand{\DJb}{\ensuremath{U_{h}^{DJ}(\borel_{+}) } }
\newcommand{\DJn}{\ensuremath{U_{h}^{DJ}(\nilp_{+}) } }

\newcommand{\kzt}{A_{\g,t}}

\newcommand{\flip}{\operatorname{\tau}}

\newcommand{\OUh}{{U}_{h}}

\newcommand{\OJ}{{J}}

\newcommand{\DJbp}{\mathcal{B}_{+}}

\newcommand{\Nform}{\mathcal{C}}
\newcommand{\E}{\mathbb{E}}
\newcommand{\F}{\mathbb{F}}

\newcommand{\DJg}{\ensuremath{U_{h}^{DJ}(\g) } }
\newcommand{\Ker}{\ensuremath{\operatorname{Ker}}}

\newtheorem{Df}{Definition}
\newtheorem{theorem}[Df]{Theorem}

\newtheorem{prop}[Df]{Proposition}

\newtheorem{lemma}[Df]{Lemma}

\newtheorem{remark}[Df]{Remark}

\newtheorem{corollary}[Df]{Corollary}
\numberwithin{Df}{section}
\numberwithin{equation}{section}

\begin{document}
\title{Some remarks on quantized Lie superalgebras of classical type}
    
\author{Nathan Geer}
\address{School of Mathematics\\ 
Georgia Institute of Technology\\ 
Atlanta, GA 30332-0160}
\email{geer@math.gatech.edu}
\date{\today}

 \begin{abstract}
In this paper we use the Etingof-Kazhdan quantization of Lie bi-superalgebras to investigate some interesting questions related to Drinfeld-Jimbo type superalgebra associated to a Lie superalgebra of classical type.   It has been shown that the D-J type superalgebra associated to a Lie superalgebra of type A-G, with the distinguished Cartan matrix, is isomorphic to the E-K quantization of the Lie superalgebra.  The first main result in the present paper is to extend this to arbitrary Cartan matrices.   This paper also contains two other main results: 1) a theorem stating that all highest weight modules of a Lie superalgebra of type A-G can be deformed to modules over the corresponding D-J type superalgebra and 2) a super version of the Drinfeld-Kohno Theorem.  
\end{abstract}

\maketitle

\section{Introduction} 
Let us start by recalling that Kac \cite{K} showed that Lie superalgebras of type A-G are characterized by their associated Dynkin diagrams or equivalently Cartan matrices.  A Cartan matrix associated to a Lie superalgebra is a pair consisting of a matrix $A$ and a set $\tau$ determining the parity of the generators.  Let $(A,\tau)$ be such a Cartan matrix and $\g$ be the Lie superalgebra arising from $(A,\tau)$.  

The Drinfeld-Jimbo algebra associated to a semi-simple Lie algebra is defined by generators and relations.  The higher order relations of this algebra are called the quantum Serre relations.  Many authors (see for example: \cite{FLV,KTol2,Yam94}) have studied generalization of this algebra to the setting of Lie superalgebras.  This generalization introduces defining relations (e.g. \eqref{E:QserreC}) that are not of the form of the quantum Serre relations.  These new relations depend directly on the Cartan matrix $(A,\tau)$.  Let $\DJgA$ be the Drinfeld-Jimbo type superalgebra associated to the triple $(\g,A,\tau)$.  

In \cite{D7} Drinfeld asked: ``Does there exist a universal quantization for Lie bialgebras?''  Etingof and Kazhdan \cite{EK1} gave a positive answer to this question.  In \cite{G04A} the author extended this quantization from the setting of Lie bialgebras to the setting of Lie bi-superalgebras.   The triple $(\g,A,\tau)$ has a natural Lie bi-superalgebra structure.  Let $\UhA$ be the Etingof-Kazhdan quantization of this Lie bi-superalgebra.  Let $\h$ be the Cartan sub-superalgebra of $\g$.  

\begin{theorem}\label{T:UhIsoDJ} (proof in \S \ref{SS:ProofThm})
There exists an isomorphism of quantized universal enveloping (QUE) superalgebras:
$$\alpha: \DJgA \rightarrow \UhA$$
such that $\alpha|_{\h}=id$.
\end{theorem}

%

In \cite{G04A} the author proves the above theorem in the case when $\g$ is a Lie superalgebra of type A-G with the distinguished Cartan matrix.  Most of the arguments given in \cite{G04A} can be adapted to the present situations.  However, in \cite{G04A} the author checks by hand that the quantum Serre relations (which depend on the Cartan matrix) are in the kernel of a certain bilinear form.  Here we appeal to the work of Yamane \cite{Yam94} to show that these relations are in the desired kernel.

  In the remainder of this section we state the other main results of this paper.
 In the case when $\g$ is a semi-simple Lie algebra, the following two theorems are analogous to results of Drinfeld \cite{D2, D5}.  Drinfeld's proof uses deformation theoretic arguments based on the fact that $H^{i}(\g,U(\g))=0, \; i=1,2,$ for semisimple Lie algebras.  In general, this vanishing result is not true for Lie superalgebras (for example $H^{1}(\mathfrak{sl}(2|1),U(\mathfrak{sl}(2|1)))\neq 0$).  Our proof is based on a different approach than Drinfeld's, utilizing Theorem \ref{T:UhIsoDJ} and the general theory of the Etingof-Kazhdan quantization of Lie (super)bialgebras (see Equation~\ref{E:DJtop}). 
 
\subsection{Deformations of weight modules} Lusztig proved that each irreducible dominant integral weight module of a Kac-Moody algebra can be deformed to a module over the corresponding Drinfeld-Jimbo algebra.  In \cite{D7} Drinfeld asked if arbitrary weight modules over a Kac-Moody algebra can be deformed.   Etingof and Kazhdan gave a positive answer to this question in \cite{EK6}.   The following theorem gives a positive answer to this question for Lie superalgebras of type A-G.

For the definition of highest weight modules over $\g$ and $\DJgA$ see subsection \ref{SS:DefWeightMod}.

\begin{theorem}\label{C:DefExist} (proof in \S\ref{SS:DefWeightMod})
For $\Lambda \in \h^{*}$, let $V(\Lambda)$ be the irreducible highest weight module over $\g$ of highest weight $\Lambda$.  Then there exists a highest weight $\DJgA$-module $\widetilde{V}(\Lambda)$ of weight $\Lambda$ which is a deformation of $V(\Lambda)$.  Moreover, the characters of $V(\Lambda)$ and $\widetilde{V}(\Lambda)$ are equal.  In other words, if
  $$V(\Lambda)   = \oplus_{\lambda \in \h^{*}}V_{\lambda} $$
  then 
  $$ \widetilde{V}(\Lambda) = \oplus_{\mu \in \h^{*}[[h]]}\widetilde{V}_{\mu}$$
  where $\widetilde{V}_{\mu}:= \{v : av=\mu(a)v \text{ for all } a\in\h\}$ and $\widetilde{V}_{\mu}\cong V_{\mu|_{\h^{*}}}[[h]]$.
\end{theorem}

Theorem \ref{C:DefExist} allows one to derive properties about $\DJgA$-modules by working with the underlying $\g$-module.  In particular, character formulas, tensor product decompositions, and other properties about $\g$-modules lead to analogous properties for the corresponding $\DJgA$-module.  This procedure is very useful in knot theory.  For example, it is used to construct generalized multivariable Alexander link invariants arising from Lie superalgebras (see \cite{GP1,GP2}).


Let use now say a few words about the proof of Theorem \ref{C:DefExist} (a detailed proof is given in Section \ref{SS:DefWeightMod}).  
In the case when $\g$ is equal to $\mathfrak{osp}(1|2n)$, $\mathfrak{osp}(2|2n)$ and $\mathfrak{sl}(m|n)$ it has been shown that $H^{2}(\g,U(\g))=0$ (see \cite{SZ,SuZ}).   Therefore, for such a Lie superalgebra $\g$ it follows that any deformation of $U(\g)$ is trivial as an associative superalgebra.  This allow $\g$-modules to be deformed to $\DJgA$-modules.  However, this does not imply that a weight $\g$-module of weight $\Lambda$ will be deformed to a weight $\DJgA$-module whose weight is equal (relative to $\h$) to $\Lambda$.   In the case of finite dimensional Lie algebras, Drinfeld gives further argument (using the vanishing of the first cohomology) to show that weight modules can be deformed (see Section 4 of \cite{D2}).  As mentioned above such arguments will not work in the present situation.

Instead our proof is based on the fact that the E-K quantization $\UhA$ is by construction the twist of a quasi-Hopf superalgebra whose underlying superalgebra is $U(\g)[[h]]$.  Combining this fact with Theorem~\ref{T:UhIsoDJ} it follows that we have an isomorphism of superalgebras
\begin{equation}
\label{E:DJtop}
f:\DJgA\rightarrow U(\g)[[h]]
\end{equation}
such that $f|_{\h}=id$.  We will show Theorem \ref{C:DefExist} follows from the observation that highest weight $\g$-modules can be deformed to highest weight $U(\g)[[h]]$-modules and the fact that $f$ preserves weights.   

\subsection{The Drinfeld-Kohno Theorem}
Here we state the Drinfeld-Kohno theorem for Lie superalgebras.  In the coming sections we elaborate on the definitions of the objects involved in this statement.

 Let $V(\Lambda)$ be an irreducible highest weight module of $\g$ and let $\widetilde{V}(\Lambda)$ be the $\DJgA$-module given in Theorem \ref{C:DefExist}. 
Let $B_{n}=\langle\sigma_{i}\rangle $ be the braid group.  
Define $\rho_{n}$ to be the representation of $B_{n}$ on $\widetilde{V}(\Lambda)^{\otimes n}$ given by
$$\sigma_{i}\mapsto \flip_{i, i+1}R_{i i+1}$$
where $\flip_{i,i+1}$ is the super permutation of the i-th and the (i+1)-th component and $R$ is the universal $R$-matrix of $\DJgA$.  Finally, let $\rho_{n}^{KZ}$ be the monodromy representation of $B_{n}$ arising from the KZ system of differential equations.

The proof of the following theorem can be found in subsection \ref{SS:D-KT}.

\begin{theorem} (The Drinfeld-Kohno theorem for Lie superalgebras) \label{C:D-KT} The representations  $\rho_{n}$ and $\rho_{n}^{KZ}$ are equivalent.
\end{theorem}


\subsection*{Acknowledgments}   I would like to thank Arkady Berenstein, Jon Brundan, Jon Kujawa and Arkady Vaintrob for helpful conversations.   I am especially grateful to Pavel Etingof.  

\section{Preliminaries}\label{Prelim}\label{S:Lbialg}

In this section we recall facts and definitions related to Lie super(bi)algebras, for more details see \cite{K}.

Let $k$ be a field of characteristic zero.  A \emph{superspace} is a $\Z_{2}$-graded vector space $V=V_{\p 0}\oplus V_{\p 1}$ over $k$.  We denote the parity of a homogeneous element $x\in V$ by $\p x\in \Z_{2}$.  We say $x$ is even (odd) if $x\in V_{\p 0}$ (resp. $x\in V_{\p 1}$).   In this paper the tensor product will have the natural induced $\Z_{2}$-grading.  Throughout, all modules will be $\Z_{2}$-graded modules, i.e. module structures which preserve the $\Z_{2}$-grading (see \cite{K}).

A \emph{Lie bi-superalgebra} is a Lie superalgebra $\g$ with a linear map $\delta : \g \rightarrow \wedge^{2}\g$ that preserves the $\Z_{2}$-grading and satisfies both the super-coJacobi identity and cocycle condition (see \cite{A}).
A triple $(\g,\g_{+},\g_{-})$ of finite dimensional Lie superalgebras is a finite dimensional \emph{super Manin triple} if $\g$ has a non-degenerate super-symmetric invariant bilinear form $\langle\, , \rangle$, such that $\g \cong \g_{+}\oplus \g_{-}$ as superspaces, and $\g_{+}$ and $\g_{-}$ are isotropic Lie sub-superalgebras of $\g$.  
There is a one-to-one correspondence between finite dimensional super Manin triples and finite dimensional Lie bi-superalgebras (see \cite[Proposition 1]{A}).


Now we give the notion of the double of a finite dimensional Lie bi-superalgebra.  Let $(\g_{+}, [\: , ]_{\g_{+}}, \delta)$ be a finite dimensional Lie bi-superalgebra and $(\g,\g_{+},\g_{-})$ its corresponding super Manin triple.  Then $\g:=\g_{+}\oplus \g_{-}$ has a natural structure of a quasitriangular Lie bi-superalgebra, see \cite{G04A}.    We call $\g$ the \emph{double} of $\g_{+}$ and denote it by $D(\g_{+})$.

We will now define the Casimir element of $\g$.  Let $p_{1},...,p_{n}$ be a homogeneous basis of $\g_{+}$.  Using the isomorphism $\g_{-} \rightarrow \g_{+}^{*}$ pick a homogeneous basis $m_{1},...,m_{n}$ of $\g_{-}$ that is dual to $p_{1},...,p_{n}$, i.e. $\langle m_{i},p_{j}\rangle=\delta_{i,j}$.  Notice $m_{1},...,m_{n},p_{1},...,p_{n}$ is a basis of $\g$ that is dual to the basis $p_{1},...,p_{n},(-1)^{\p m_{1}}m_{1},...,(-1)^{\p m_{n}}m_{n}$, with respect to $\langle\, , \rangle$.  Define the \emph{Casimir element} to be 
\begin{equation}
\label{D:Cas}
\cas =\sum p_{i}\otimes m_{i}+ \sum (-1)^{\p m_{i}}m_{i}\otimes p_{i}.
\end{equation}  The element $\cas$ is even, invariant and super-symmetric.  Moreover, it is independent of the choice of basis.

\section{The Drinfeld-Jimbo type quantization of Lie superalgebras of type A-G}\label{S:QUg}

In this section we recall some basic facts related to complex Lie superalgebras of type A-G and their quantum analogue. 
For the purposes of this paper Lie superalgebras of type A-G will be complex and include the Lie superalgebra $D(2,1,\alpha)$.

Any two Borel subalgebras of a semisimple Lie algebra are conjugate.  Moreover, semisimple Lie algebras are determined by their root systems or equivalently their Dynkin diagrams.   Not all Borel sub-superalgebras of classical Lie superalgebras are conjugate.  As shown by Kac \cite{K} a Lie superalgebra can have more than one Dynkin diagram depending on the choice of Borel.  However, using Dynkin type diagrams Kac gave a characterzation of Lie superalgebras of type A-G.  The constructions of this paper depend on the choice of Borel sub-superalgebra.   

Let $h$ be an indeterminate. 

\subsection{Lie superalgebras of type A-G}\label{SS:g}   

Let $\g$ be a Lie superalgebra of type A-G.  Let $\Phi=\{\alpha_{1},...,\alpha_{s}\}$ be a simple root system of $\g$ and let $(A,\tau)$ be the corresponding Cartan matrix, where $A$ is a $s\times s$ matrix and $\tau$ is a subset of $\{1,...,s\}$ determining the parity of the generators.  The matrix $A=(a_{ij})$ is symmetrizable, i.e. there exists nonzero rational numbers $d_{1},\dots,d_{s}$ such that $d_{i}a_{ij}=d_{j}a_{ji}$.  By rescaling, if necessary, we may and will assume that $d_{1}=1$.  For notational convenience we set $I=\{1,...,s\}$.  

From Propositions 2.5.3 and 2.5.5 of \cite{K} there exists a unique (up to constant factor) non-degenerate supersymmetric invariant bilinear form $(,)$ on $\g$.  Moreover, the restriction of this form to the Cartan sub-superalgebra $\h$ is non-degenerate.  Let $h_{i}$, $i\in I$, be defined by $(a,h_{i})=d_{i}^{-1}\alpha_{i}(a)$ for all $a\in \h$.

Yamane \cite{Yam94} showed that $\g$ is given by generators and relations which depend on $(A,\tau)$.  We will now recall this presentation.   

The Lie superalgebra $\g$ is generated by $h_{i}$, $e_{i},$ and $f_{i}$ for $i\in I$ (whose parities are all even except for  ${e}_{i}$ and ${f}_{i}$ if $i\in \tau$ which are odd) such that the relations satisfy: 
\begin{align}
\label{R:LieSuperalg}
   [h_{i}, h_{j}]&=0,   & [h_{i}, e_{j}] &=a_{ij}e_{j},   &
   [ h_{i},f_{j}] &=-a_{ij} f_{j} &   [e_{i},f_{j}] &=\delta_{ij}h_{i}
\end{align}
and the \emph{super Serre} relations
$$[e_{i},e_{i}] =[f_{i},f_{i}] =0 \text{ for } i\in \tau$$
$$(\text{ad } e_{i})^{1+|a_{ij}|}e_{j}=(\text{ad } f_{i})^{1+|a_{ij}|}f_{j}=0, \text{ if } i \neq j,\text{ and } i \notin \tau $$
plus nonstandard super Serre-type relations which depend on $(A,\tau)$.  For example 
\begin{equation*}
[e_{m},[e_{m-1},[e_{m},e_{m+1}]]]=[f_{m},[f_{m-1},[f_{m},f_{m+1}]]]=0
\end{equation*}
if $ m-1,m, m+1\in I $, $  a_{mm}=0$ and the Cartan matrix $A$ is not of type C or D.  
For a complete list of relations see \cite{Yam94}.

We denote a Lie superalgebra of type A-G by the triple $(\g,A,\tau)$.


\subsection{Lie bi-superalgebra structure}\label{SS:glbialg} 
In this subsection we will recall that Lie superalgebras of type A-G has a natural Lie bi-superalgebra structure.  Similar results have previously been considered by other authors (see for example \cite{BGN,LS}).  The results of this subsection are straightforward generalizations of the non-super case.  

Let $(\g,A,\tau)$ be a  Lie superalgebra of type A-G.  Let $\h=\langle h_{i}\rangle _{i \in I}$ be the Cartan subalgebra of $\g$.
Let $\nilp_{+} \: (\text{resp., } \nilp_{-}) $ be the nilpotent Lie sub-superalgebra of $\g$ generated by $e_{i} \text{'s }$ $(\text{resp., } f_{i}\text{'s})$.  Let $\borel_{\pm} := \nilp_{\pm}\oplus \h $ be the Borel Lie sub-superalgebra of $\g$.  

Let $\eta_{\pm}:\borel_{\pm} \rightarrow \g \oplus \h$ be defined by
$$\eta_{\pm}(x)=x\oplus (\pm \bar x),$$
where $\bar x$ is the image of $x$ in $\h$.  Using this embedding we can regard $\borel_{+} $ and $\borel_{-}$ as Lie sub-superalgebras of $\g \oplus \h$


As above let $(,)$ be the unique non-degenerate supersymmetric invariant bilinear form on $\g$. 
 Let $(,)_{\g \oplus \h}:=(,)-(,)_{\h}$, where $(,)_{\h}$ is the restriction of $(,)$ to $\h$.


\begin{prop}\label{P:glManinT}
$(\g \oplus \h, \borel_{+}, \borel_{-})$ is a super Manin triple with $(,)_{\g \oplus \h}$.  
\end{prop}

\begin{proof}
Under the embedding $\eta_{\pm}$ the Lie subsuperalgebra $\borel_{\pm}$ is isotropic with respect to $(,)_{\g \oplus \h}$.   Since $(,) $ and $(,)_{\h}$ both are invariant super-symmetric nondegenerate bilinear forms then so is $(,)_{\g \oplus \h}$.  Therefore the Proposition follows.
\end{proof}

The Proposition implies that $\g \oplus \h, \borel_{+}$ and $ \borel_{-}$ are Lie bi-superalgebras.  Moreover, we have that $\borel_{+}^{*}\cong\borel_{-}^{op}$ as Lie bi-superalgebras, where $^{op}$ is the opposite cobracket.  

We will now compute the formulas for the cobrackets of these Lie bi-superalgebras.  The cobracket on $\borel_{+}\subset \g\oplus \h$ is induced by $\borel_{-}$ under the pairing $(,)_{\g \oplus \h}:\borel_{+}\otimes \borel_{-}\rightarrow \C$.  In other words, this pairing induces a Lie superalgebra structure on $\borel_{+}^{*}$ and thus a Lie supercoalgebra structure on $\borel_{+}^{**}\cong \borel_{+}$.  Using these facts we compute the cobracket on $\borel_{+}$.  

Let $f_{i}$ and $K_{i}$ be a basis of $\borel_{-}$ so that  
$$(h_{i}\oplus h_{i},K_{j})_{\g \oplus \h}=\delta_{ij}$$
where $h_{i}\oplus h_{i} \in \h \subset \borel_{+}$.  Let $\pi:\g\oplus \h\rightarrow \g$ be the natural projection.  Set $k_{j}:=\pi(K_{j})$.  Then 
$$(h_{i}\oplus h_{i},K_{j})_{\g \oplus \h}=2(h_{i},k_{j})=2d_{i}^{-1}\alpha_{i}(k_{j}),$$
implying $\alpha_{i}(k_{j})=\delta_{ij}d_{i}/2$.  By definition
$$[h^{*}_{j},e^{*}_{i}]:= [k_{j},f_{i}]=-\alpha_{i}(k_{j})f_{i}=-\delta_{ij}\frac{d_{i}}{2}e^{*}_{i}.$$  
From the discussion above we see that 
\begin{equation}
\label{E:cobracketb+}
\delta(e_{i})=\frac{d_{i}}{2}(e_{i}\otimes h_{i}-h_{i}\otimes e_{i})=\frac{d_{i}}{2}e_{i}\wedge h_{i}.
\end{equation}
Similarly,
\begin{align}
\label{E:cobracketb-}
\delta(f_{i})&=-\frac{d_{i}}{2}f_{i}\wedge h_{i} & \delta(a) &=0 \text{ for } a \in \h\subset \borel_{\pm}
\end{align}
These formulas define a Lie bi-superalgebra structure on the Lie superalgebras $\g, \borel_{+}$ and $\borel_{-}$.  Moreover, we have that ($\g,\bar{r})$ is a quasitriangular Lie bi-superalgebra where $\bar{r}$ is the image of the canonical element $r$ in the double $D(\borel_{+})\cong \g \oplus \h$ under the natural projection.


\subsection{The Drinfeld-Jimbo type superalgebra $\DJgA$} \label{SS:GenRel}  
Khoroshkin-Tolstoy \cite{KTol} and Yamane \cite{Yam94} used the quantum double notion to define a quasitriangular QUE superalgebra $\DJgA$.  In this subsection we recall their results which are needed in this paper.

Set $q=e^{h/2}$ and $q_{i}=q^{d_{i}}$.  


\begin{theorem}[\cite{KTol,Yam94}] 
Let $(\g,A,\tau)$ be a Lie superalgebra of type A-G.  There exists an explicit quasitriangular QUE superalgebra $(\DJgA, R)$.  The superalgebra $\DJgA$ is defined as the $\C[[h]]$-superalgebra generated by $\h$ and the elements $\E_{i}$ and $\F_{i}, $ $i \in I$ (all generators are even except $\E_{i}$ and $\F_{i}$ for $i\in \tau$ which are odd) 
\begin{align}
\label{E:DJglRelation1}
 [a,a']  &=0, & [a,\E_{i}]=&\alpha_{i}(a)\E_{i}, & [a,\F_{i}]=&-\alpha_{i}(a)\F_{i}, \text{ for } a, a' \in \h  
\end{align}
\begin{align}
\label{E:DJglRelation2}
 [\E_{i},\F_{j}]=&\delta_{i,j}\frac{q^{d_{i}h_{i}}-q^{-d_{i}h_{i}}}{q_{i}-q_{i}^{-1}},  
 \end{align}
plus the super quantum Serre-type relations (see \cite{Yam94}).
 The coproduct and counit given by
\begin{align*}
\label{}
    \D {\E_{i}}= & \E_{i}\otimes q^{d_{i}h_{i}} + 1 \otimes \E_{i}, & \D {\F_{i}}  = & \F_{i}\otimes 1 + q^{-d_{i}h_{i}} \otimes \F_{i},  \\
   \D a = & a \otimes 1 + 1\otimes a, &  \epsilon(a) = &  \epsilon(\E_{i})= \epsilon(\F_{i})=0
\end{align*}
for all $ a\in \h$.  Moreover, the $R$-matrix is given by explicit formulas.
\end{theorem}
We call $\DJgA$ the Drinfeld-Jimbo type quantization of $(\g,A,\tau)$.
\begin{remark}
The super quantum Serre-type relations depend directly on $(A,\tau)$.  For example, 
\begin{multline}
\label{E:QserreC}
\E_{m}\E_{m-1}\E_{m}\E_{m+1}+\E_{m}\E_{m+1}\E_{m}\E_{m-1}+\E_{m-1}\E_{m}\E_{m+1}\E_{m}
\\+\E_{m+1}\E_{m}\E_{m-1}\E_{m}-(q+q^{-1})\E_{m}\E_{m-1}\E_{m+1}\E_{m}=0
\end{multline}
if $ m-1,m, m+1\in I $, $  a_{mm}=0$ and the Cartan matrix $A$ is not of type C or D.  
\end{remark}
\begin{remark}\label{R:GenRelinKer}
In defining $\DJgA$, Yamane constructed a bilinear form on a free algebra whose kernel is the so called super quantum Serre-type relations.  We will now give the properties of this form, which we will use later.  

Let $\DJbp$ 
be the Hopf superalgebra generated over $\C[[h]]$ by $\h$ and $E_{i}$, $i \in I$, where $E_{i}$ is odd if $i\in \tau$ and all other generators are even, with relations satisfying
\begin{align*}
\label{}
    [a,a']=&0, & [a,E_{i}]=&\alpha_{i}(a)E_{i},
\end{align*}  
and coproduct defined by
\begin{align*}
\label{}
    \D {E_{i}}= & E_{i}\otimes q^{d_{i}h_{i}} + 1 \otimes E_{i} &
   \D a = & a \otimes 1 + 1\otimes a
\end{align*}
for $ a, a'\in \h$ and $i\in I$.  In \cite{Yam94} Yamane defined a $\C((h))$-valued bilinear form $\Nform$ on $\DJbp$ with the following properties:
\begin{align*}
\label{}
  \Nform(xy,z)=& \Nform(x \otimes y,\D z), & \Nform(x, yz)=& \Nform(\D x, y \otimes z) 
\end{align*}
(where $\Nform(x\otimes y,z\otimes w):=(-1)^{\p y \p z}\Nform(x,z)\Nform(y,w)$) and
$$ \Nform(E_{i},E_{j})=
	\begin{cases}
		(q_{i}-q_{i}^{-1})(q^{d_{\alpha}}-q^{-d_{\alpha}})^{-1} & \text{if $i=j$,}\\
		0 & \text{otherwise}
	\end{cases}$$
	where $d_{\alpha}=1$ if $a_{ii}=0$ and $d_{\alpha}=d_{i}a_{ii}/2$ otherwise.  
Yamane showed that $\Ker(\Nform)$ is generated by the super quantum Serre-type relations (see sections 4.2 and 10.4 of \cite{Yam94}).  Moreover, the morphism $\DJbp\rightarrow \DJgA$ which is the identity on $\h$ and maps $E_{i}$ to $\E_{i}$ induces an isomorphism $\DJbp/\Ker(\Nform)\rightarrow \DJb$, where $\DJb$ is the sub-superalgebra of $\DJgA$ generated by $\h$ and $\E_{i}$, $i \in I$.  
\end{remark}


\section{The Etingof-Kazhdan quantization}
In this section, we will show that for Lie superalgebras of type A-G, the Etingof-Kazhdan quantization is isomorphic to the Drinfeld-Jimbo quantization.  As in \cite{EK6} we will show that the E-K quantization is given by the desired generators and relations.  In particular, we show that the E-K quantization of the Borel sub-bi-superalgebra is isomorphic to a superalgebra given by generators and relations modulo the kernel of an appropriate bilinear form.    It then follows from Yamane's work that the E-K quantization is in fact the D-J quantization.  Throughout this section we will use the notation of section \ref{S:QUg}.


Let $\ghat$ be the Lie superalgebra generated by $e_{i}, f_{i}$ and $h_{i}$ for $i\in I$ satisfying relation (\ref{R:LieSuperalg}) where all generators are even except ${e}_{i}$ and ${f}_{i}$ for $i\in \tau$ which are odd.  Let $\borelhat_{+}$ and $\borelhat_{-}$ be the Borel sub-superalgebras of $\ghat$ generated by $e_{i}, h_{i}$ and $f_{i}, h_{i}$, respectively.  
The formulas \eqref{E:cobracketb+} and \eqref{E:cobracketb-} define Lie bi-superalgebra structures on $\ghat$ and $\borelhat_{\pm}$. 

Let $\Uh$ be the Etingof-Kazhdan quantization of a Lie bi-superalgebra $\g$ defined in \cite{G04A}.  This quantization has two important properties that we will use here.  The first is that the quantization is functorial.  The second is that it commutes with taking the double, i.e. $$D(\Uh)\cong U_{h}(D(\g))$$
where $D(\Uh)$ is the quantum double and $D(\g)$ is the double of $\g$ (see \S \ref{S:Lbialg}).  When $\g$ is a Lie superalgebra of type A-G its natural bi-superalgebra structure depends on the choice of Cartan matrix (see \S\ref{SS:glbialg}).  For this reason we denote the E-K quantization of such an Lie superalgebra by $\UhA$.


\subsection{Generators and relations for $\Uhbhat$}

\begin{theorem}\label{T:GRofUhbhat}
The quantized universal enveloping (QUE) superalgebra $\Uhbhat$ is isomorphic to the QUE superalgebra $\Uhatp$ generated over $\C[[h]]$ by $\h$ and the elements $E_{i}$ for $i\in I$ (all generators are even except for $E_{i}$, $i\in \tau$ which are odd) satisfying the relations
\begin{align*}
\label{}
    [a,a']=&0, & [a,E_{i}]=&\alpha_{i}(a)E_{i},
\end{align*}  
with coproduct
\begin{align*}
\label{}
    \D{a}=&1\otimes a + a \otimes 1, & \D{E_{i}}=&E_{i}\otimes q^{d_{i}h_{i}}+ 1\otimes E_{i},
\end{align*}   
for all $a,a'\in \h$ and $i,j\in I$.
\end{theorem}

The theorem follows from the following two lemmas.  


\begin{lemma}\label{L:UQESbhat}
The QUE superalgebra $\Uhbhat$ is isomorphic to the QUE superalgebra generated over $\C[[h]]$ by $\h$ and the elements $E_{i}, \: i\in I$ (all generators are even except for $E_{i}$, $i\in \tau$ which are odd) satisfying the relations
\begin{align}
\label{E:relUhbhatCar}
[a,a']=& 0
\end{align}
\begin{align}
\label{E:relUhbhat}
       [a,E_{i}]=&\alpha_{i}(a)E_{i},
\end{align}  
with coproduct
\begin{align}
\label{E:coprdCar}
 \D{a}=&1\otimes a + a \otimes 1, 
\end{align}
\begin{align}
\label{E:coprd}
  \D{E_{i}}=&E_{i}\otimes q^{\gamma_{i}}+ 1\otimes E_{i},
\end{align}   
for all $a,a'\in \h$ and $i,j\in I$ and suitable elements $\gamma_{i}\in \h[[h]]$.
\end{lemma}

\begin{proof}
Since the E-K quantization is functorial, the embedding of Lie bi-superalgebras $\h \rightarrow \borelhat_{+}$ induces a embedding of QUE superalgebras $\Uhh \rightarrow \Uhbhat$.  Note that this embedding of QUE superalgebras restricted to $\h$ is the identity.  We will use this observation later. 

 By construction $\Uhh$ is equal to $U(\h)[[h]]$.  It follows that $\Uhbhat$ has a sub-superalgebra generated by $\h$ which satisfies relation \eqref{E:relUhbhatCar} and whose coproduct is given by \eqref{E:coprdCar}.
Since $\borelhat_{+}=\h\oplus \tilde{\nilp}_{+}$ where $ \tilde{\nilp}_{+}$ is the free Lie bi-superalgebra generated by $e_{i}$, $i\in I$ and $\Uhbhat\cong U(\borelhat)[[h]]$ as a superalgebra, to complete the proof, it suffices to show that there exists $E_{i}$ in $\Uhbhat$ which satisfy relations \eqref{E:relUhbhat}  and \eqref{E:coprd}.   

The Lie bi-superalgebra $\borelhat_{+}$ has a natural $\Z_{+}^{n}$-grading given by $\mathrm{deg}_{i}(h_{j})=0$ and $\mathrm{deg}_{i}(e_{j})=\delta_{ij}$.  The functorality of the quantization implies $\Uhbhat$ has a $\Z_{+}^{n}$-grading.  It follows that $\Uhbhat=\oplus_{\textbf{m}\in \Z_{+}^{n}}\Uhbhat [\textbf{m}]$
where $\Uhbhat [\textbf{m}]$ is a free $\Uhh$-module of finite rank.  In particular, let $1_{i}\in\Z_{+}^{n}$ be given by $1_{i}(j)=\delta_{ij}$ then $\Uhbhat [1_{i}]$ has rank 1.
Let $E_{i}'$ be an element in $\Uhbhat[1_{i}]$ such that $E_{i}'$ is $e_{i}$ modulo $h$.

The proof is completed by showing that there exist an element $x$ in $1+hU(\h)[[h]]\subset \Uhh$ such that $E_{i}=E_{i}'x$ satisfies  \eqref{E:coprd}.  After replacing the ordinary tensor product with the super-tensor product, the construction of $x$ follows as in the proof of Proposition 3.1 of \cite{EK6}.  There are no new signs introduced.  
For the most part, this is true because the arguments of the proof are based on the purely even Cartan subalgebra $\h$ and the functorality of the quantization.  
\end{proof}


\begin{lemma}\label{L:gamma}
$ \gamma_{i}=d_{i}h_{i}$
\end{lemma}
\begin{proof} 
By definition we have the natural projection $\borelhat_{+}\rightarrow\borel_{+}$.  Then the functoriality of the quantization implies that there is an epimorphism of Hopf superalgebras $U_{h}(\borelhat_{+}) \rightarrow U_{h}(\borel_{+})$.  Therefore, $U_{h}(\borel_{+})$ is generated by $\h$ and $E_{i}$ satisfying the relations \ref{E:relUhbhatCar}-\ref{E:coprd} (and possibly other relations).  So it suffices to show that $ \gamma_{i}=d_{i}h_{i}$ in $U_{h}(\borel_{+})$.

Next we show that $U_{h}(\borel_{+}) \cong U_{-h}(\borel_{+})\qdo$ where ${}\qdo$ denotes the quantum dual with the opposite coproduct (see \cite{G04A}).  From the definition of $\g$ the Lie bi-superalgebra $\borel_{+}$ is self dual, i.e. $\borel_{+}\cong \borel_{+}^{*}$.  Again from functoriality we have that $U_{h}(\borel_{+}) \cong U_{h}(\borel_{+}^{*})$.  From Proposition \ref{P:glManinT} we have $\borel_{+}^{*}\cong\borel_{-}^{op}$.  Then equation (45) of \cite{G04A} imply that $U_{h}(\borel_{+})\qdo \cong U_{h}(\borel_{+}^{*op})$.  Substituting $\borel_{+}^{op}$ for $\borel_{+}$ we have $U_{h}(\borel_{+}^{op})\qdo \cong U_{h}(\borel_{+}^{*})$.  Finally from relation (7) of \cite{G04A} it follows that $U_{h}(\borel_{+}^{op}) \cong U_{-h}(\borel_{+})$ which implies that $U_{-h}(\borel_{+})\qdo \cong U_{h}(\borel_{+}^{*})$.  Thus, we have shown that $U_{h}(\borel_{+}) \cong U_{-h}(\borel_{+})\qdo$.  

This isomorphism gives rise to the bilinear form $B: U_{h}(\borelp) \otimes U_{-h}(\borelp) \rightarrow \C((h))$ which satisfies the following conditions 
\begin{align}
\label{E:RelationsOfB1}
  B(xy,z)=& B(x \otimes y,\D z), & B(x, yz)=& B(\D x, y \otimes z) 
\end{align}
$$ B(q^{a},q^{b})=q^{-(a,b)}, a,b \in \h.$$
Let $a\in\h$ and $i\in I$.  Set $B_{i}=B(E_{i},E_{i})$, which is nonzero.   Using (\ref{E:RelationsOfB1}) we have 
\begin{align*}
\label{}
    B(E_{i},q^{a}E_{i})&=B(E_{i}\otimes q^{ \gamma_{i}} +1 \otimes E_{i}, q^{a}\otimes E_{i})   \\
    &  = B(E_{i}, q^{a})B(q^{ \gamma_{i}},E_{i}) + B(1,q^{a})B(E_{i},E_{i}) \\
    & = B_{i}
\end{align*}
since $B(E_{i},q^{a})=0$.  Similarly, we have $B(E_{i},q^{a}E_{i}q^{-a})= B(E_{i},q^{a}E_{i})B(q^{\gamma_{i}},q^{-a})$ implying
\begin{align}
\label{E:GammaB}
    B_{i}q^{(a,\gamma_{i})}=&B(E_{i},q^{a}E_{i}q^{-a}).  
\end{align}  
To complete the proof we need the following relation:
\begin{align}
\label{E:quantea}
   q^{a}E_{i}q^{-a}=&q^{\alpha_{i}(a)}E_{i} 
\end{align}
This relation is equivalent to $q^{h_{j}}E_{i}q^{-h_{j}}=q^{\alpha_{i}(h_{j})}E_{i}$ which follows from expanding $q=e^{h}$ and using the relation $[a,E_{i}]=\alpha_{i}(a)E_{i}$.
From (\ref{E:GammaB}) and (\ref{E:quantea}) we have 
$$B_{i}q^{(a,\gamma_{i})}=B(E_{i},q^{a}E_{i}q^{-a})=B(E_{i},q^{\alpha_{i}(a)}E_{i})=B_{i}q^{\alpha_{i}(a)}.$$
Thus, $(a, \gamma_{i})=\alpha_{i}(a)$, but $ \alpha_{i}(a)= d_{i}(a,h_{i})$, and so $\gamma_{i}=d_{i}h_{i}$, which completes the proof. 
\end{proof}


\subsection{The quantized universal enveloping superalgebra $\Uhb$}

In this subsection we show that there exist a bilinear form on $\Uhbhat$ such that $\Uhbhat$ modulo the kernel of the form is isomorphic to $\Uhb$.  


\begin{theorem}\label{T:FormB}
There exists a unique bilinear form on $\Uhbhat$ which takes values in $\C((h))$ with the following properties
\begin{align*}
\label{}
  B(xy,z)=& B(x \otimes y,\D z), & B(x, yz)=& B(\D x, y \otimes z) 
\end{align*}
$$ B(q^{a},q^{b})=q^{-(a,b)}, a,b \in \h,$$
$$ B(E_{i},E_{j})=
	\begin{cases}
		(q_{i}-q_{i}^{-1})(q^{d_{\alpha}}-q^{-d_{\alpha}})^{-1} & \text{if $i=j$,}\\
		0 & \text{otherwise}
	\end{cases}$$
	where $d_{\alpha}=1$ if $a_{ii}=0$ and $d_{\alpha}=d_{i}a_{ii}/2$ otherwise.  Moreover $\Uhb \cong \Up:=\Uhatp/\Ker(B)$ as QUE superalgebras.	
\end{theorem}
\begin{proof}
The existence and uniqueness follows from the fact that the superalgebra generated by the $E_{i}$ is free. 

We will show that there is a nondegenerate bilinear form on $\Uhb$ with the same properties as $B$.  From the proof of Lemma \ref{L:gamma} we have that $U_{h}(\borel_{+}) \cong U_{-h}(\borel_{+})\qdo$.  But the even homomorphism $U_{-h}(\borel_{+})^{op} \rightarrow \Uhb$ given by conjugation by $q^{-\sum x_{i}^{2}/2}$, where $x_{i}$ is a orthonormal basis for $\h$, is a isomorphism.  Therefore we have a even isomorphism $\Uhb \cong \Uhb\qd$.  This isomorphism gives rise to the desired form on $\Uhb$.  

So the form $B$ is the pull back of the form on $\Uhb$.  Implying that the kernel of the form on $\Uhb$ is contained in the image of the kernel of $B$ under natural projection.

But the kernel of the form on $\Uhb$ is zero since the form is nondegenerate.   Thus we have $\Uhb \cong \Uhatp/\Ker(B)$. 
\end{proof}

\begin{corollary}\label{C:QUEBorelIso}
Let $\DJb$ be the sub-superalgebra of $\DJgA$ generated by $\h$ and the elements $\E_{i}$, $i \in I$.  Then the map 
$$g: \DJb\rightarrow \Uhb \text{ given by } g|_{\h}=id  \text{ and } \E_{i}\mapsto E_{i}$$ 
is an isomorphism of QUE superalgebras. 
\end{corollary}
\begin{proof}
The corollary follows directly from Theorem \ref{T:FormB} and Remark \ref{R:GenRelinKer}.
\end{proof}


\subsection{Proof of Theorem \ref{T:UhIsoDJ}}\label{SS:ProofThm}

\begin{proof}[Proof of Theorem \ref{T:UhIsoDJ}]
It suffices to show that the QUE superalgebra $\UhA$ is isomorphic to the quotient of the double $D(\Up)$ by the ideal generated by the identification of $\h \subset \Up$ and $\h^{*} \subset \Upqd$.  

Recall from \S \ref{SS:glbialg} that the Lie bi-superalgebra structure of $\g$ comes from identifying $\h$ and $\h^{*}$ in $\g \oplus \h= \borel_{+}\oplus \borel_{+}^{*}$.  Also since the quantization commutes with the double we have 
$$\OUh(D(\borel_{+})) \cong D(\Uhb)=\Uhb \otimes \Uhb\qdo.$$ 
Therefore, we have $\UhA$ is isomorphic to $D(\Uhb)=\Uhb \otimes \Uhb\qdo$ modulo the the ideal generated by the identification of $\h \subset \Uhb$ and $\h^{*} \subset \Uhb\qdo$.  But from Theorem \ref{T:FormB} we have that $D(\Uhb)\cong D(\Up)$ and then Corollary \ref{C:QUEBorelIso} implies result.
\end{proof}

\subsection{Twisting of Drinfeld-Jimbo superalgebras}
In this subsection we give a corollary of Theorem \ref{T:UhIsoDJ}. 

Khoroshkin and Tolstoy \cite{KTol2} showed that any two isomorphic Lie superalgebras with different Cartan matrices have isomorphic deformations (as associative superalgebras) and their coproducts are connected by a twisting of a factor of the universal $R$-matrix.  It is not clear if the definition of the quantized superalgebra associated to a Lie superalgebra of type A-G is correct in \cite{KTol2} (as some relations appear to be missing).  However, after adding the missing relations it becomes apparent that the results of \cite{KTol2} hold.  

In any case, we give an alternative proof that any two isomorphic Lie superalgebras with different Cartan matrices have isomorphic Drinfeld-Jimbo type superalgebras (as associative superalgebras) and their coproducts are connected by a twist.  The primary difference in our approach, as opposed to \cite{KTol2}, is to construct the twist using the E-K quantization rather than as a factor of the universal $R$-matrix. 

\begin{corollary}\label{C:IsoForDifBorel}
Let $\g$ and $\g'$ be two isomorphic Lie superalgebras of type A-G with associated Cartan matrices $(A,\tau)$ and $(A',\tau')$.  Let $\OJ$ and $\OJ'$ be the element of $\U[[h]]^{\otimes 2}$ defined in (32) of \cite{G04A} using $(\g,A,\tau)$ and $(\g',A',\tau')$, respectively.  Then the QUE superalgebra $\DJgA$ is isomorphic to $\DJgAp$ twisted by the element $(\OJ')^{-1}\OJ$.  In particular, $\DJgA$ and $\DJgAp$ are isomorphic as associative superalgebras.   
\end{corollary}
\begin{proof}
Let $\cas$ and $\cas'$ be the Casimir elements of $(\g,A,\tau)$ and $(\g',A',\tau')$, respectively.  Since the Casimir element is independent of the choice of basis we have $A_{\g,\cas}$ and $A_{\g',\cas'}$ are isomorphic quasitriangular quasi-Hopf superalgebras.  Recall that by construction $\UhAp=(A_{\g',\cas'})_{\OJ'}$, and so 
$$ \left(\UhAp\right)_{\OJ'^{-1}}=A_{\g',\cas'}.$$
  Similarly, $\UhA=(A_{\g,\cas})_{\OJ}$ implying 
 \begin{equation}
\label{E:uhtwist}
 \UhA\cong\left((\UhAp)_{\OJ'^{-1}}\right)_{\OJ} 
\end{equation}
  as quasitriangular quasi-Hopf superalgebras.  But $\UhA$ is a quasitriangular Hopf superalgebra and so the result follows from Equation~\eqref{E:uhtwist} and Theorem~\ref{T:UhIsoDJ}.  
\end{proof}

As mentioned above the relations of the D-J type superalgebra depend on the choice of the Cartan matrix.  For this reason it is not apparent from the definition that $\DJgA$ and $\DJgAp$ are isomorphic as associative superalgebras.

\section{A theorem of Drinfeld's}\label{S:ModuleCat} Let $(\g,A,\tau)$ be a Lie superalgebra of type A-G.  Recall from section \ref{S:Lbialg} that for each super Manin triple there exists a Casimir element.  Let $\cas$ be this element associated to the triple $(\g,\borel_{+},\borel_{-})$.  

For each Lie algebra and symmetric invariant 2-tensor Drinfeld \cite{D3} defined a quasitriangular quasi-Hopf quantized universal enveloping algebra:  
$$(U(\g)[[h]],\Delta_{0},\epsilon_{0},R_{KZ}=e^{th/2},\Phi_{KZ}).$$
The morphisms $\Delta_{0}$ and $\epsilon_{0}$ are the standard coproduct and counit of $U(\g)[[h]]$.  The element $\Phi_{KZ}$ is the KZ-associator.  
Setting $t=\cas$ let $\kzt$ be the analogous topologically free quasitriangular quasi-Hopf superalgebra (for more details see \cite{G04A,Z}).  Also recall the definition $\DJgA$ given in \S\ref{SS:GenRel}.  

Here we show that the categories of topologically free modules over $\kzt$ and $\DJgA$ are braided tensor equivalent.  We do this in two steps: (1) we show that $\UhA$ and $\kzt$ have equivalent module categories, (2) we use the fact the that $\DJgA$ and $\UhA$ are isomorphic to prove the desired result.   For more on braided tensor categories see \cite{Kas}. 

\subsection{The E-K quantization $\UhA$ and $\kzt$}
In this subsection we show that $\UhA\text{-}Mod$  and $\kzt\text{-}Mod$ are equivalent braided tensor categories.  To this end we recall the following definitions.  

Let $(A,\Delta,\epsilon, \Phi,R)$ be a quasitriangular quasi-bi-superalgebra (see \S4.1 of \cite{G04A}).  An invertible element $J\in A\otimes A$ is a \emph{gauge transformation on $A$} if 
$$(\epsilon \otimes id)(J)=(id \otimes \epsilon)(J)=1.$$
Using a gauge transformation $J$ on $A$, one can construct a new quasitriangular quasi-bi-superalgebra $A_{J}$ with coproduct $\Delta_{J}$, R-matrix $R_{J}$ and associator $\Phi_{J}$ defined by
$$\Delta_{J}= J^{-1}\Delta J, \:\: R_{J}=(J^{op})^{-1}RJ,$$
$$  \Phi_{J}= J^{-1}_{23}(id \otimes \Delta)(J^{-1})\Phi(\Delta \otimes id)(J)J_{12}.$$

As is the case of quasitriangular (quasi-)bialgebra, the category of modules over a quasitriangular (quasi-)bi-superalgebra is a braided tensor category.  Let $X$ be a topological (quasi-)bi-superalgebra and let $ X \text{-}Mod$ be the category of topologically free $X$-modules. 

\begin{theorem}\label{T:EquivOfqqbialg}
Let $A$ and $A'$ be a quasitriangular quasi-bi-superalgebra.  Suppose that $J$ is a gauge transformation on $A'$ and $\alpha:A\rightarrow A_{J}'$ is an isomorphism of quasitriangular quasi-bi-superalgebras.  Then $\alpha$ induces a equivalence between the braided tensor categories $A'\text{-}Mod$ and $A\text{-}Mod.$
\end{theorem}
\begin{proof}
Let $\alpha^{*}:A'\text{-}Mod \rightarrow A\text{-}Mod$ be the functor defined as follows.  On objects, the functor $\alpha^{*}$ is defined by sending the module $W$ to the same underlying vector space with the action given via the isomorphism $\alpha$.  For any morphism $f: W\rightarrow X$ in $A'\text{-}Mod$ let $\alpha^{*}(f)$ be the image of $f$ under the isomorphism
$$\Hom_{A'}(W,X)\cong \Hom_{A}(W,X).$$
A standard categorical argument shows that this functor is an equivalence of braided tensor categories (see \S XV.3 of \cite{Kas}).
\end{proof}

Let $\OJ$ be the element of $\U[[h]]^{\otimes 2}$ defined in (32) of \cite{G04A}.  The definition of the element $\OJ$ uses the associator $\Phi$.  By construction the E-K quantization is the twist of $\kzt$ by the element $\OJ$, i.e. $\UhA=(\kzt)_{\OJ}$.  For exact formulas of the coproduct and R-matrix of $\UhA$ see Proposition 16 and the end of section 5.2 of \cite{G04A}.

\subsection{A braided tensor equivalence}\label{SS:BTE}

The following theorem was first due to Drinfeld \cite{D5} in the case of semi-simple Lie algebras.

\begin{theorem}\label{T:Drinfeldsuper}
The braided tensor categories $\kzt\text{-}Mod$ and $\DJgA\text{-}Mod$ are equivalent. 
\end{theorem}

\begin{proof}
As mentioned at the end of the last subsection $\UhA=(\kzt)_{\OJ}$.  Combining this fact with 
Theorem \ref{T:UhIsoDJ} we have that there exists an isomorphism of quasitriangular quasi-bi-superalgebra
$$\alpha: \DJgA \rightarrow (\kzt)_{\OJ}.$$ 
Now as a consequence of Theorem \ref{T:EquivOfqqbialg} we have that $\kzt\text{-}Mod$ and $\DJg\text{-}Mod$ are braided tensor equivalent.
\end{proof}

\begin{remark}
Drinfeld's proof of Theorem \ref{T:Drinfeldsuper} in the case of semi-simple Lie algebras uses deformation theoretic arguments to show the existence of $\alpha$.  Our proof constructs the isomorphism $\alpha$ explicitly.
\end{remark}

\section{Proofs} We will now give the proofs of Theorems \ref{C:DefExist} and \ref{C:D-KT}. 

\subsection{Proof of Theorem \ref{C:DefExist}}\label{SS:DefWeightMod}  In this subsection we give the definitions of highest weight modules and a proof of Theorem \ref{C:DefExist}.

Let $(\g,A,\tau)$ be a Lie superalgebra of type A-G and let $\h$ and $\nilp_{+}$ be its Cartan and nilpotent sub-superalgebras, respectively.  Let $\Lambda$ be an element of $ h^{*}$.  Let $c_{\Lambda}$ be the one dimensional $\borel_{+}$-module generated by $v_{\Lambda}$ with the following action:
  \begin{align*}
\label{}
   \nilp_{+}v_{\Lambda} &=0, & av_{\Lambda}=\Lambda(a)v_{\Lambda} \text{ for } a \in \h
\end{align*}
and where we set $\p{v}_{\Lambda}=\p 0$.  Set $\widehat{V}(\Lambda)=\mathrm{Ind}_{\borel_{+}}^{\g}c_{\Lambda}$.  Then $\widehat{V}(\Lambda)$ contains a unique maximal submodule $I(\Lambda)$.  We call $V(\Lambda):=\widehat{V}(\Lambda)/I(\Lambda)$ the irreducible weight module with highest weight $\Lambda$.

Next we define a similar notion for $\DJgA$-module.  A topologically free $\DJgA$-module $V$ is call a highest weight module with highest weight $\Lambda\in \h^{*}$ if there exists a non-zero even generating vector $v_{\Lambda}$ such that 
\begin{align*}
\label{}
  \DJn v_{\Lambda} &=0, & av_{\Lambda}=\Lambda(a)v_{\Lambda} \text{ for } a \in \h    
\end{align*}
where $\DJn$ is the sub-superalgebra of $\DJgA$ generated by $\E_{i}$, $i \in I$. 

\begin{proof}[Proof of Theorem \ref{C:DefExist}]
Consider the $U(\g)[[h]]$-module $V(\Lambda)[[h]]$.  From Theorem \ref{T:Drinfeldsuper} we know that $\alpha$ induces a $\DJgA$-module structure on $V(\Lambda)[[h]]$, denote this $\DJgA$-module by $\widetilde{V}(\Lambda)$.  In addition, from Corollary \ref{C:QUEBorelIso} we have $\alpha|_{\h}=Id$ and $\alpha$ restricted to $U(\borel_{+})[[h]]$ is an isomorphism between $U(\borel_{+})[[h]]$ and $\DJb$.  Thus, $\widetilde{V}(\Lambda)$ is a highest weight $\DJgA$-module such that $\widetilde{V}(\Lambda)/h\widetilde{V}(\Lambda)\cong V(\Lambda)$.
\end{proof}


\subsection{Proof of Theorem \ref{C:D-KT}}\label{SS:D-KT} Here we give the proof of the Drinfeld-Kohno Theorem for Lie superalgebras.  Before proving the theorem we will define the $KZ$ monodromy representation of the braid group.

Let $(\g,A,\tau)$ be a Lie superalgebra of type A-G.  Let $V$ be a irreducible highest weight module of $\g$.  Consider the Knizhnik-Zamolodchikov system of differential equations with respect to a function $\omega(z_{1},...,z_{n}) $ of complex variables $z_{1},...,z_{n}$ with values in $V^{\otimes n}[[h]]$:
\begin{equation}
\label{E:KZvalueV}
\frac{\partial \omega}{\partial z_{i}} = \frac{\hbar}{2\pi i}\sum_{i\neq j}\frac{\cas_{ij}\omega}{z_{i}-z_{j}}.
\end{equation}

We have that this system of equations defines a flat connection on the trivial bundle $Y_{n}\times V^{\otimes n}[[h]]$
where $Y_n=\{(z_1,...,z_n) | i\neq j \text{ implies } z_i \neq z_j \}\subset \C^{n}.$
This connection determines a monodromy representation from $\pi_{1}(Y_{n})$ to $\Aut_{\C[[h]]}(V^{\otimes n}[[h]])$.  Moreover, since the system of equations \eqref{E:KZvalueV} is invariant under the action of the symmetric group we obtain a monodromy representation 
$$\rho_{n}^{KZ}: \pi_{1}(X_{n},p)\rightarrow \Aut_{\C[[h]]}(V^{\otimes n}[[h]])$$
where $X_{n}=Y_{n}/S_{n}$ and $p=(1,2,...,n)\in \C^{n}$.  Finally, we identify $\pi_{1}(X_{n},p)$ with the braid group $B_{n}$ to get a monodromy representation of $B_{n}$.

\begin{proof}[Proof of Theorem \ref{C:D-KT}]
Let $\rho_{n}^{R_{KZ}}$ be the representation of $B_{n}$ on $V(\Lambda)^{\otimes n}[[h]]$ induced by the $R$-matrix $R_{KZ}=e^{h\cas/2}$.  From Theorem \ref{T:Drinfeldsuper} and Theorem \ref{C:DefExist} we have that $\rho_{n}^{R_{KZ}}$ and $\rho_{n}$ correspond to each other under the braided tensor functor $\alpha^{*}$.  The theorem follows since $\rho^{KZ}_{n}$ coincides with $\rho_{n}^{R_{KZ}}$. 
\end{proof}

\end{document}